\documentclass[12pt]{article}    

\input{amssym.def} 
\input{amssym.tex}

\begin{document}      
%**********************************************    
\title{Hyperbolic geometry, continued fractions
      and classification of the finitely generated totally ordered
      simple dimension groups}    
%***********************************************    

\author{Igor  ~V.  ~Nikolaev
\footnote{Partially supported 
by NSERC.}}

\date{}     
%\date{March 20, 2003}    
 \maketitle

\newtheorem{thm}{Theorem}    
\newtheorem{lem}{Lemma}    
\newtheorem{dfn}{Definition}    
\newtheorem{rmk}{Remark}    
\newtheorem{cor}{Corollary}    
\newtheorem{prp}{Proposition}    
\newtheorem{exm}{Example}    
    
%******************************************************************    
 \begin{abstract}    
We  classify the polycyclic totally ordered simple dimension groups,  i.e.  dimension groups 
given  by  a dense embedding of  lattice  ${\Bbb Z}^{n}$ into  the real line.  Our method is based on
the  geometry of  simple geodesics  on  the hyperbolic surface  of genus $g\ge2$.
The main theorem says that isomorphism classes of the  polycyclic totally ordered dimension groups 
are bijective   with a generic subset of  reals $\alpha$  modulo the action of  group  $GL(2,{\Bbb Z})$. 
The result is an extension of the Effros-Shen classification of the dicyclic dimension
groups.

\vspace{7mm}    
    
{\it Key words and phrases:  dimension group,  geodesic lamination,  Jacobi-Perron continued fraction}

\vspace{5mm}
{\it MSC.:  46L85 (noncommutative topology); 57M50  (geometric structures on low-dimensional manifolds)}
\end{abstract}

%**************************************************************************
\section{Introduction}
%***************************************************************************
The dimension of a subspace of the Euclidean space ${\Bbb R}^n$  can be 1, 2, 3 or more.
It is no longer true that  dimension of any non-trivial subspace of the infinite-dimensional  Euclidean
space ${\Bbb R}^{\infty}$ is a positive integer.  It was discovered by von Neumann 
that there exist linear subspaces in ${\Bbb R}^{\infty}$ of a non-integer 
(continuous) dimension.   Namely,  there exists a dimension function on the projections
in a von Neumann algebra  ranging  in the unit interval.

Unlike the von Neumann algebras,   the  dimension function on the  projections of a  
$C^*$-algebra, $A$,   takes value in an abelian group, $K_0(A)$,  rather than in ${\Bbb R}$,
see e.g.  [Effros  1981]   \cite{E},   Chapter  1.   For instance,  if $A$  is an AF-algebra {\it ibid.},    then
the range of the dimension function in $K_0(A)$ is known as a {\it dimension
group} of $A$,  see appendix for the definition.  The classification of such groups is a difficult open problem.

In  the remarkable paper [Effros  \&  Shen  1980]   \cite{EfS},   the authors  classified  the dicyclic dimension groups,
i.e.  dimension groups inside the  abelian group  ${\Bbb Z}^2$. Their main  result says that each simple 
dicyclic group can be assigned a  positive irrational number $\alpha$,  defined  modulo  the action
of matrix group $GL(2,{\Bbb Z})$,  and such that if
$\alpha=[a_0,a_1,a_2,\dots]$
is a regular continued fraction of $\alpha$,  then one gets a representation   
of the dicyclic group via the simplicial dimension groups:
%**************************************************************
\begin{equation}\label{eq2}
{\Bbb Z}^2\buildrel\rm 
\left(\small\matrix{0 & 1\cr 1 & a_0}\right)
\over\longrightarrow {\Bbb Z}^2
\buildrel\rm 
\left(\small\matrix{0 & 1\cr 1 & a_1}\right)
\over\longrightarrow
{\Bbb Z}^2\buildrel\rm 
\left(\small\matrix{0 & 1\cr 1 & a_2}\right)
\over\longrightarrow \dots
\end{equation}
%*************************************************************** 
see [Effros  \&  Shen  1980]   \cite{EfS},   Theorem 3.2.
The  irrational $\alpha$ is a {\it slope}  of the  straight line in the plane, which is
the universal cover of   two-dimensional torus, see [Effros  \&  Shen  1980]   \cite{EfS} Theorem 2.1. 
The range of the dimension function (a positive cone) in ${\Bbb Z}^2$ is the half-plane 
$P_{\alpha}=\{(k,l)\in {\Bbb Z}^2~|~\alpha k+l\ge 0\}$,  and it is shown that each  simple
dicyclic dimension group  arises in this way.

The objective of  our  note is    similar classification of the polycyclic
dimension groups,  i.e. dimension groups inside the abelian group ${\Bbb Z}^n$.   
Recall that the  universal cover for  surfaces of genus $g\ge 2$ is no longer the Euclidean plane
but the hyperbolic (Lobachevsky) half-plane  ${\Bbb H}$.  It is well known,  that the role of 
straight lines in ${\Bbb H}$ is played  by the simple geodesics, i.e. 
geodesics with no self-crossing points;  we refer the reader to appendix for an exact definition.   
Thus,  to classify the polycyclic dimension groups,  one needs: 

\medskip
(i) to define a dimension  group, $G$,  coming from the simple
geodesic, $\gamma$,  on a hyperbolic surface of genus $g\ge2$;

\smallskip
(ii) to define a slope $\alpha\in {\Bbb R}$ of $\gamma$ on the surface;

\smallskip 
 (iii) to construct  a simplicial   approximation   of dimension group $G$ 
in terms of the  slope $\alpha$ of geodesic $\gamma$ .

\bigskip\noindent
The realization of (i) -- (iii) is as follows. The closure,  $\bar\gamma$,
of a simple non-periodic geodesic $\gamma$ consists of the continuum 
of disjoint non-periodic geodesics known as a {\it geodesic lamination}
$\lambda$ [Casson \& Bleiler  1988]   \cite{CaB}.  Let $|\lambda|$ be the total number of the principal
regions of $\lambda$ ({\it ibid.}, p.60) and $n=2g+|\lambda|-1$, where
$g$ is the genus of the hyperbolic surface, $X$,  carrying the lamination $\lambda$. 
It is  known,  that the set of invariant transversal measures of $\lambda$
is a convex compact set  $\Delta_{k-1}$  of dimension $1\le k\le\left[{n\over 2}\right]$,
where $[\bullet]$ is the integer part of a number, see  [Sataev  1975]    \cite{Sat}.  By $G_{\lambda}$
we understand the dimension group inside the abelian group ${\Bbb Z}^n$,
whose state space $S(G_{\lambda})$ is isomorphic to $\Delta_{k-1}$. 
The assignment of $G_{\lambda}$ via the invariant measures on $\lambda$ is unique [Goodearl  1986]    \cite{G}, Chapter 4. 
To implement (ii),  we shall use the  {\it Jacobi-Perron continued fraction}  attached 
to the geodesic lamination $\lambda$.  Namely, a standard bijection (a blow-up, see
the appendix) between geodesic laminations and foliations gives
rise to  a measured foliation ${\cal F}_{\lambda}$ on the surface $X$. 
 For simplicity, we let ${\cal F}_{\lambda}$
be given by trajectories of a closed 1-form  $\omega_{\lambda}$  on $X$. 
Let $\lambda_i=\int_{\gamma_i}\omega_{\lambda}$ be the periods  of
 $\omega_{\lambda}$ in a basis $\{\gamma_1,\dots,\gamma_n\}$
of the relative homology  $H_1(X, Sing~\omega_{\lambda}; {\Bbb Z})$;
choosing  $\lambda_1\ne 0$,   we let $\theta_{i-1}=\lambda_i /\lambda_1$ for   $i\ge 2$ 
and $\theta=(\theta_1,\dots,\theta_{n-1})$.   Consider  the Jacobi-Perron continued 
fraction of  $\theta$ [Bernstein  1971]   \cite{B}:
%**********************************************************************************
\begin{equation}\label{eq3}
\left(\matrix{1\cr \theta}\right)=
\lim_{k\to\infty} \left(\matrix{0 & 1\cr I & b_0}\right)\dots
\left(\matrix{0 & 1\cr I & b_k}\right)
\left(\matrix{0\cr {\Bbb I}}\right),
\end{equation}
%************************************************************ 
where $I$ is the unit matrix, $b_i=(b^{(i)}_1,\dots, b^{(i)}_{n-1})^T$ a 
vector of  the non-negative integers and ${\Bbb I}=(0,\dots, 0, 1)^T$.
By a {\it slope} of the geodesic $\gamma$, we shall understand 
an irrational number  $\alpha$  given by the regular continued fraction:
%**********************************************************************************
\begin{equation}\label{eq3bis}
\alpha=[b_1^{(0)}+1,\dots, b_{n-1}^{(0)}+1, b_1^{(1)}+1,\dots, b_{n-1}^{(1)}+1,\dots].
\end{equation}
%************************************************************ 
(In other words, $\alpha$ is a concatenation from zero to infinity of the vectors $b_i$
with the  entries added by 1 and,  thus,  all  strictly positive.) 
Finally, to realize item (iii),  one needs
to restrict to the generic case $k=1$,  i.e. the totally ordered dimension
groups.  
Such a restriction  secures convergence of the Jacobi-Perron
continued fraction (\ref{eq3}) 
and one can apply the known theorem of  [Effros  \&  Shen  1979]   \cite{EfSh2} on 
approximation of  the unimodular dimension groups by the simplicial 
dimension groups.  A summary of our results can be formulated as follows.
%*****************************************************************
\begin{thm}\label{thm1}
Each finitely generated totally ordered simple dimension group $G$
can be indexed by a positive irrational number $\alpha\in U$,
where $U$ is a generic subset of ${\Bbb R}$. The parametrization, $G_{\alpha}$,
has the following properties:

\medskip
(i) if $G_{\alpha}$ and $G_{\alpha'}$ are order-isomorphic,  then:
%*****************************************************************************

\begin{equation}
\alpha'={a\alpha+b\over c\alpha+d}  \quad\hbox{for a matrix} 
\quad \left(\matrix{a & b\cr c & d}\right) \in GL(2, {\Bbb Z});
\end{equation}
%**************************************************************************

\smallskip
(ii) $G_{\alpha}$  is the limit of the following  simplicial dimension groups: 
%*************************************************************
\begin{equation}\label{eq4}
{\Bbb Z}^n
\buildrel\rm
\left(\matrix{0 & 1\cr I & b_0}\right)
\over\longrightarrow 
{\Bbb Z}^n
\buildrel\rm
\left(\matrix{0 & 1\cr I & b_1}\right)
\over\longrightarrow
{\Bbb Z}^n
\dots
\end{equation}
%**************************************************************
\end{thm}
%***********************************************************
The article is organized as follows.
 The dimension group $G_{\lambda}$  is introduced in section 2. 
Theorem \ref{thm1}  is proved in section 3. Finally, in section 4 a brief account 
of geodesic laminations,  measured foliations and the Jacobi-Perron fractions is 
given;  it is supplemented by basic facts on the dimension groups, whose  
 complete coverage can be found in  [R\o rdam, Larsen \& Laustsen 2000]    \cite{RLL}, Chapter 5.

%**********************************************************************
\section{Dimension groups generated by foliations}
%***********************************************************************
For notation, we refer the reader to appendix;  
let $\gamma$ be a simple non-periodic geodesic on a hyperbolic  surface $X$
of genus $g\ge 2$. The closure $\bar\gamma$ contains a continuum of non-periodic
pairwise disjoint simple geodeisics, which form a perfect (Cantor) set on  $X$. 
Such a closure is a geodesic lamination, which we shall denote
by $\lambda$. The number of principal regions of $\lambda$ will be denoted 
by $|\lambda|$.

It is well known that $\lambda$ corresponds to a foliation ${\cal F}_{\lambda}$
on $X$, see e.g. [Thurston  1997]    \cite{T}, Chapter 8.5. 
The foliation ${\cal F}_{\lambda}$ has  $|\lambda|$ singular
points of  saddle type; it  is obtained  from the geodesic lamination  $\lambda$
by a blow-down homotopy (see the appendix). 
Denote by $\Phi_X$ the space of foliations on $X$, whose singularity set coincides
with that  of ${\cal F}_{\lambda}$. The coordinates $(\lambda_1,\dots,\lambda_n)$ of 
${\cal F}_{\lambda}$ in $\Phi_X$ are given by the formula: 
%**************************************************************
\begin{equation}\label{eq6}
\lambda_i=\int_{\gamma_i} \omega_{\lambda} ~d\mu,
\end{equation}
%*************************************************************
where $\omega_{\lambda}$ is a closed 1-form tangent to the leaves of
foliation ${\cal F}_{\lambda}$, 
$\{\gamma_1,\dots,\gamma_n\}$ is a basis in the relative homology
$H_1(X, Sing~{\cal F}_{\lambda}; {\Bbb Z})$ and $\mu$ 
an invariant transversal measure on the leaves of foliation
${\cal F}_{\lambda}$, see Section 4.2 for the details.  It follows from the formulas
for the relative homology that:  
%**********************************************************
\begin{equation}\label{eq7}
n=2g+|\lambda|-1.
\end{equation}
%************************************************************
On the other hand,   it is known that the total number of independent
invariant measures of foliation ${\cal F}_{\lambda}$ is equal to
$k\le \left[{n\over 2}\right]$,  see  [Sataev  1975]   \cite{Sat}.
If we denote by $(\mu_1,\dots,\mu_k)$ all such measures,  then
parameter $\mu$  in formula  (\ref{eq6}) can be written as 
 $\mu=\sum_{j=1}^k\alpha_j\mu_j$,
where $\alpha_j\ge 0$ are some real numbers. 
Thus,  one gets from formula  (\ref{eq6})  an $n\times k$ matrix $(\lambda_{ij})$, whose entry
$\lambda_{ij}$ is obtained by integration of $\omega$
along the contour $\gamma_i$ with  $\mu=\mu_j$;  each $j$-row 
of the matrix defines a homomorphism $h_j: {\Bbb Z}^n\to {\Bbb R}$,
such that $h_j(e_i)=\lambda_{ij}$, where $e_i$ is a base element in 
${\Bbb Z}^n$.  It is easy to see,  that the kernel of $h_j$
is a hyperplane in ${\Bbb R}^n$. Thus, one obtains a dimension group, $G_{\lambda}$,
inside the abelian group ${\Bbb Z}^n$, which is bounded by the
$k$ hyperplanes corresponding to the measures $(\mu_1,\dots,\mu_k)$.    
%*************************************************************************
\begin{dfn}\label{dfn1}
The dimension group $G_{\lambda}$ is called associated to the
geodesic lamination $\lambda$ (equivalently, foliation ${\cal F}_{\lambda}$). 
\end{dfn}
%***********************************************************************
In the sequel, our main case will be $k=1$ (i.e. the totally ordered dimension group
$G_{\lambda}$). Note that the foliations ${\cal F}_{\lambda}$ with a unique invariant ergodic
measure are generic in the space $\Phi_X$ [Masur  1982]    \cite{Mas1}.  In this generic case, 
the dimension group $G_{\lambda}$ can be  identified with the projective class 
of a  ${\Bbb Z}$-module  ${\Bbb Z}\lambda_1+\dots+{\Bbb Z}\lambda_n\subset {\Bbb R}$,
i.e. an equivalence class of the modules $\mu({\Bbb Z}\lambda_1+\dots+{\Bbb Z}\lambda_n)$,
where $\mu>0$ is a real number.

%**************************************************************************
\section{Proof of theorem \ref{thm1}}
%***************************************************************************
To  prove theorem \ref{thm1},  we have to show that   parametrization $G_{\alpha}$
of the polycyclic  totally ordered simple dimension groups described
in introduction  has the following  properties:

\medskip
(i) if $G_{\alpha}$ and $G_{\alpha'}$ are order-isomorphic,  then
%*****************************************************************************

\displaymath
\alpha'={a\alpha+b\over c\alpha+d}  \quad\hbox{for a matrix} 
\quad \left(\matrix{a & b\cr c & d}\right) \in GL(2, {\Bbb Z});
\enddisplaymath
%**************************************************************************

\smallskip
(ii) $G_{\alpha}$  is the limit of the following  simplicial dimension groups: 
%*************************************************************
\displaymath
{\Bbb Z}^n
\buildrel\rm
\left(\matrix{0 & 1\cr I & b_0}\right)
\over\longrightarrow 
{\Bbb Z}^n
\buildrel\rm
\left(\matrix{0 & 1\cr I & b_1}\right)
\over\longrightarrow
{\Bbb Z}^n
\dots
\enddisplaymath
%**************************************************************

\bigskip
We shall proceed stepwise;  we refer the reader to the appendix for all preliminary facts and notation
used below.

\medskip
(i) Let $G$ be a finitely generated totally ordered simple dimension group.
In this case $G\cong {\Bbb Z}^n$ is a dense subgroup 
of  the real line ${\Bbb R}$.  Let $\lambda_1,\dots,\lambda_n$ be a set of
generators of the subgroup;  one  can always choose $\lambda_i$ to be positive.

We wish to construct a measured foliation ${\cal F}$  from the set of reals 
$\lambda_1,\dots,\lambda_n$. The foliation is carried by a surface $X$, whose
genus we shall further specify. Our goal can be achieved with the help 
of the method of zippered rectangles (see the appendix).
First, we define the bottom  of the rectangle to be $\lambda_1+\dots+\lambda_n$;
so far the construction is not unique pending a choice of a permutation $\pi$.
Next,  one takes an irreducible permutation $\pi$ on the $n$ symbols;  the 
top of the rectangle will be  $\lambda_{\pi(1)}+\dots+\lambda_{\pi(n)}$.
(If one takes a different $\pi$, one would get another parameterization 
of dimension groups by measured foliations;  note, however, that the number
of such choices is always finite.)   
Instead of $\pi$, one can specify the number and types of the singular
points of the foliation ${\cal F}$; there exists a one-to-one
correspondence between these two sets of data [Masur  1982]    \cite{Mas1}, [Veech  1982]    \cite{Vee1}.
Now, the method of zippered rectangles  produces a unique measured foliation ${\cal F}$ 
on a  surface $X$, whose genus is given by the formula: 
%**************************************************************
\begin{equation}
g={1\over 2}(n+1-|Sing~{\cal F}|),
\end{equation}
%***************************************************************
where $|Sing~{\cal F}|$ is equal to the number of cyclic 
permutations in the prime decomposition of $\pi$.

Using the integration formula (\ref{eq17}),
it is immediate that  dimension group,  generated by the foliation ${\cal F}$,  
is order-isomorphic to $G$.  Note that since
$G$ is totally ordered, the foliation ${\cal F}$ is uniquely ergodic.

We can now index the group $G$ by a slope $\alpha$ of the foliation $\cal F$;
we shall write the corresponding dimension group as $G_{\alpha}$.   
Denote by $i$ an (index) map which assigns to every totally ordered simple dimension
group a slope $\alpha$ of the corresponding measured foliation ${\cal F}$.

Let us prove the following series of lemmas, which reflect the properties of the
index map. 
%******************************************************
\begin{lem}\label{lm1}
The index map $i$ is an injection.  
\end{lem}
%******************************************************
{\it Proof.}
Indeed, let to the contrary, $G\ne G'$ and $i(G)=i(G')$. 
By the uniqueness theorem ([Perron  1907]    \cite{Per1}, Section 5, Theorem 4), the Jacobi-Perron 
fractions for the vectors $\lambda, \lambda'$ must be different
in at least one term. So does the regular continued fractions,
which define the slopes $\alpha,\alpha'$. Thus, by a main property
of the continued  fractions, we have $\alpha\ne\alpha'$. One runs into
a contradiction with the initial assumption.
$\square$

%******************************************************
\begin{lem}\label{lm2}
If the totally ordered dimension groups $G_{\alpha}$ and $G_{\beta}$ are
order-isomorphic,  then 
%******************************************************
\begin{equation}\label{eq12}
\beta={a\alpha+b\over c\alpha +d},
\end{equation}
%******************************************************                      
where $a,b,c,d$ are integers such that $ad-bc=\pm 1$. 
\end{lem}
%*****************************************************
{\it Proof.}
Let $G_{\alpha}$ and $G_{\beta}$ be given by the ${\Bbb Z}$-modules
$\sum {\Bbb Z}\lambda_i$ and $\sum {\Bbb Z}\lambda_i'$, respectively.
Since $G_{\alpha}$ and $G_{\beta}$ are order-isomorphic, after a scaling:
%*********************************************************
\begin{equation}
{\Bbb Z}\lambda_1+\dots+{\Bbb Z}\lambda_n= {\Bbb Z}\lambda_1'+\dots+{\Bbb Z}\lambda_n',
\end{equation}
%******************************************************** 
as the subsets of the real line;   moreover,  there exists a positive isomorphism
$G_{\alpha}\to G_{\beta}$ given by the formula $\lambda_j'=\sum_{i=1}^na_{ij}\lambda_i$,
where $A=(a_{ij})$ is  invertible matrix with the non-negative integer entries $a_{ij}\ge 0$,  see e.g.
[Effros 1981]  \cite{E}, p.10.

According to [Bauer  1996]   \cite{Bau1}, Proposition 3,  matrix $A$ can be uniquely 
factorized as:
%**********************************************************************************
\begin{equation}
A=
\left(\matrix{0 & 1\cr I & b_0'}\right)\dots
\left(\matrix{0 & 1\cr I & b_s'}\right)
:= B_1'\dots B_s',
\end{equation}
%************************************************************ 
where $b_i'=(\tilde b_1^{(i)},\dots, \tilde b_{n-1}^{(i)})$ are vectors of the 
non-negative integers. Thus, $\lambda'=B_1'\dots B_s'\lambda$ and by the definition of a
slope:
%**********************************************************************************
\begin{equation}
\beta=[\tilde b^{(0)}_1 +1,\dots, \tilde b^{(s)}_{n-1} +1, \alpha]={a\alpha+b\over c\alpha +d},
\end{equation}
%************************************************************ 
for a matrix $\left(\small\matrix{a & b\cr c & d}\right)\in GL(2, {\Bbb Z})$.
$\square$

%*********************************************************************************
\begin{rmk}
\textnormal{
Note that the converse of lemma \ref{lm2} is true only in the case 
of the dicyclic dimension groups. The effect is due to the fact that the Jacobi-Perron
algorithm converges only on a generic subset of the slopes $\alpha\in {\Bbb R}$. 
}
\end{rmk}
%***********************************************************************************

\smallskip
The item (i) of theorem \ref{thm1} is proved.

\bigskip\noindent
(ii) Let ${\Bbb Z}\lambda_1+\dots+{\Bbb Z}\lambda_n$ be a ${\Bbb Z}$-module
in ${\Bbb R}$, corresponding to the totally ordered dimension group $G_{\alpha}$.
%***************************************************************************************
\begin{rmk}
\textnormal{
The  $\lambda_i$ can be obtained  from the continued fraction of $\alpha$ using an 
inverse of formula (\ref{eq3bis}) and  the Jacobi-Perron continued fraction:
%**********************************************************************************
\begin{equation}\label{eq16}
\left(\matrix{1\cr \theta}\right)=
\lim_{k\to\infty} \left(\matrix{0 & 1\cr I & b_0}\right)\dots
\left(\matrix{0 & 1\cr I & b_k}\right)
\left(\matrix{0\cr {\Bbb I}}\right),
\end{equation}
%*********************************************************************************** 
where $\theta=(\lambda_2/\lambda_1,\dots,\lambda_n/\lambda_1)$ and
$\lambda_1\ne0$.
 It is known [Bauer  1996]   \cite{Bau1}, that the fraction is convergent
on a generic set $V\subset {\Bbb R}^n$ of vectors $\lambda$ corresponding 
to the foliation ${\cal F}_{\alpha}$ with the unique ergodic measure.
In particular, the totally ordered dimension groups are in  bijection  with
such foliations and therefore our fraction is always convergent.
}
\end{rmk}
%**********************************************************************************
One can now apply a result of Effros and Shen: the dimension group $G_{\alpha}$
is a unimodular dimension group given by the following limit of simplicial dimension 
groups:
%*************************************************************
\begin{equation}\label{eq17'}
{\Bbb Z}^n
\buildrel\rm
\left(\matrix{0 & 1\cr I & b_0}\right)
\over\longrightarrow 
{\Bbb Z}^n
\buildrel\rm
\left(\matrix{0 & 1\cr I & b_1}\right)
\over\longrightarrow
{\Bbb Z}^n
\dots
\end{equation}
%**************************************************************  
see [Effros  \&  Shen  1979]   \cite{EfSh2}, Corollary 3.3. The item (ii) follows.
Theorem \ref{thm1} is proved. 
$\square$

%**************************************************************************
\section{Appendix}
%***************************************************************************
The present section contains a brief account on the geodesic laminations,
measured foliations, method of zippered rectangles  and the Jacobi-Perron continued fractions. The corresponding 
topics are covered in [Casson \& Bleiler  1988]   \cite{CaB}, [Masur  1982]    \cite{Mas1} and [Bernstein  1971]   \cite{B},  respectively. 
We also add a paragraph on the dimension groups; we refer the reader to [R\o rdam, Larsen \& Laustsen 2000]    \cite{RLL} for a 
complete account.

%**************************************************************************
\subsection{Geodesics laminations}
%***************************************************************************
Let $S$ be a finite genus surface of  constant negative curvature (a hyperbolic surface).
By a {\it geodesic}, one understands  a  maximal arc on $S$ consisting of
the locally shortest sub-arcs (in the given metric on $S$). 
Each geodesic is the image of the open real interval $I$ under a continuous map $I\to S$.
%**************************************************************************
\begin{lem}\label{lm3}
{\bf (Topological classification of geodesics)}
Let $p\in S$ be a point with an attached unit vector $t\in S^1$ on the surface
$S$. Then :

\medskip
(a) for almost all points $t\in S^1$ (w.r.t. to the Lebesgue measure on $S^1$),  the geodesic 
line through $p$ in direction $t$ is an immersion $I\to S$, i.e. a finite or infinite curve
with  self-intersections;

\smallskip
(b) the remaining  set $K\subset S^1$ has the cardinality
of continuum and geodesic lines through $p$ in the direction $t\in K$
are embedded curves of one of the three types:

\medskip\hskip1cm
(i) periodic;

\smallskip\hskip1cm
(ii) spiraling towards a periodic geodesic;

\smallskip\hskip1cm
(iii) non-periodic,  whose closure is a perfect (Cantor) subset of $S$.
\end{lem}
%**************************************************************************
{\it Proof.} See [Artin 1924]   \cite{Art} and [Myrberg  1931]    \cite{Myr}.
$\square$

\smallskip
The  geodesics of type (b) are called {\it simple}, since they have
no self-crossing points. 
Every geodesic $\gamma: I\to S$ of type (iii) is recurrent, i.e.
for any $t_0\in I$ and $\varepsilon>0$ the $\varepsilon$-neighbourhood
of $p(t_0)$ has infinitely many intersections with $\gamma(t)$ provided
$t>N$, where $N=N(\varepsilon)$ is sufficiently large. 
The topological closure of recurrent geodesic on $S$ contains
 a continuum of the disjoint recurrent geodesics,
and called a {\it geodesic lamination} $\lambda$.
The intersection of $\lambda$ with any closed curve on $S$ is a Cantor set. 
The set $S - \lambda$ is called 
a {\it principal region} of $\lambda$. The principal region can have up to 
$4g-4$ connected components on the surface of genus $g\ge2$ [Casson \& Bleiler  1988]   \cite{CaB}. 
We denote by $|\lambda|$ the total number of such components.

A {\it foliation} $\cal F$ on a surface $X$ is a partition of $X$ into a disjoint union
of $1$-dimensional and, possibly, a finite number of $0$-dimensional
leaves denoted by $Sing~{\cal F}$. The immediate examples of foliations
are orbits of the flows and trajectories of the quadratic differentials $f(z)dz^2$
on $X$ [Fathi,  Laudenbach  \&  Po\'enaru  1979]    \cite{FLP}.  The foliation $\cal F$ is called {\it measured} if it supports an invariant transversal
measure on the leaves [Fathi,  Laudenbach  \&  Po\'enaru  1979]    \cite{FLP}. In other words:

\medskip
(i) $Sing~{\cal F}$ consists of the $n$-prong saddles, where $n\ge 3$;

\medskip
(ii) each $1$-leaf is everywhere dense in $X$.

\smallskip\noindent
The geodesic lamination $\lambda$ can be obtained from $\cal F$ by a {\it blow-up}
homotopy. Namely, a {\it separatrix} of $\cal F$ is a $1$-leaf one of whose
ends lie in $Sing~{\cal F}$. The blow-up is a replacement of the separatrix
by a narrow strip $[-\varepsilon,\varepsilon]\times {\Bbb R}$ using a 
homotopy surgery, which does not affect the nearby leaves. The complement
of the blown-up separatrices consists of leaves of $\cal F$ that
make up a perfect (Cantor) set on  $X$. 
It is not hard to prove (e.g. [Thurston  1997]    \cite{T}),  that the above complement
is homeomorphic to $\lambda$.   
Moreover, $|\lambda|$ is equal to the number of singular points
of the foliation ${\cal F}$.  Note that each measured foliation
can be given by the orbits of  a closed one-form $\omega$ on $X$,
passing, if necessary,  to a double cover of $X$ [Masur  1982]    \cite{Mas1}.

%**************************************************************************
\subsection{Method of zippered rectangles}
%***************************************************************************
There exists a remarkable construction, which allows to 
produce a measured foliation from  a given set of positive reals $(\lambda_1,\dots,\lambda_n)$.  
Let $\pi$ be a permutation of $n$ symbols. Consider a rectangle with the 
base $\lambda_1+\dots+\lambda_n$ and the top $\lambda_{\pi(1)}+\dots+\lambda_{\pi(n)}$.
We shall identify the open interval $(\lambda_{i-1},\lambda_i)$ in the base with the open interval
$(\lambda_{\pi(i)-1}, \lambda_{\pi(i)})$ at the top for all $i=1,\dots, n$.  
The resulting object will be a $k$-holed topological surface, $X$, of genus
$g={1\over 2}(n-N(\pi)+1)$, where $N(\pi)$ is the number of cyclic 
permutations in the prime decomposition of $\pi$ [Veech  1982]    \cite{Vee1}.
A foliation ${\cal F}$ on $X$ is defined by vertical lines given by
the closed 1-form $\omega=dx$. The order of the singular points of ${\cal F}$
depends on the length of the elementary cyclic permutations and the total
number of the singular points equals $k=N(\pi)$. The singular points
are located at the holes of surface $X$.  
To recover $\lambda_i$ from the 1-form $\omega$, notice that
%****************************************************
\begin{equation}
n=2g+N(\pi)-1=dim~H_1(X, Sing~{\cal F}; {\Bbb Z}),
\end{equation}
%***************************************************
where the last symbol stays for the relative homology of $X$
with respect to the set of singular points of  ${\cal F}$. 
Since  $\omega=dx$, one arrives at the elementary, but important
formula:
%********************************************************
\begin{equation}\label{eq17}
\lambda_i=\int_{\gamma_i}\omega,
\end{equation}
%********************************************************
where $\gamma_i$ are the elements of a basis in $H_1(X, Sing~{\cal F}; {\Bbb Z})$.

%**************************************************************************
\subsection{Jacobi-Perron fractions}
%***************************************************************************
The Jacobi-Perron algorithm and connected (multidimensional) continued 
fraction generalizes the euclidean algorithm (regular continued fraction)
of an irrational number. Namely, let $\lambda=(\lambda_1,\dots,\lambda_n)$,
$\lambda_i\in {\Bbb R}-{\Bbb Q}$ and  $\theta_{i-1}={\lambda_i\over\lambda_1}$, where
$1\le i\le n$.   The continued fraction 
%**************************************************************
$$
\left(\matrix{1\cr \theta_1\cr\vdots\cr\theta_{n-1}} \right)=
\lim_{k\to\infty} 
\left(\matrix{0 &  0 & \dots & 0 & 1\cr
              1 &  0 & \dots & 0 & b_1^{(1)}\cr
              \vdots &\vdots & &\vdots &\vdots\cr
              0 &  0 & \dots & 1 & b_{n-1}^{(1)}}\right)
\dots 
\left(\matrix{0 &  0 & \dots & 0 & 1\cr
              1 &  0 & \dots & 0 & b_1^{(k)}\cr
              \vdots &\vdots & &\vdots &\vdots\cr
              0 &  0 & \dots & 1 & b_{n-1}^{(k)}}\right)
\left(\matrix{0\cr 0\cr\vdots\cr 1} \right),
$$
%**************************************************************
where $b_i^{(j)}\in {\Bbb N}\cup\{0\}$, is called a {\it  Jacobi-Perron
fraction},  see [Perron 1907]  \cite{Per1}. To recover the integers $b_i^{(k)}$ from the vector $(\theta_1,\dots,\theta_{n-1})$,
one has to repeatedly solve the following system of equations:
%*******************************************************************
\begin{equation}
\left\{
\begin{array}{ccc}
\theta_1 &= b_1^{(1)} &+  {1\over\theta_{n-1}^{\prime}} \nonumber\\
\theta_2 &= b_2^{(1)} &+  {\theta_1'\over\theta_{n-1}^{\prime}}\nonumber\\
\vdots & \nonumber\\
\theta_{n-1} &= b_{n-1}^{(1)} &+  {\theta_{n-2}'\over\theta_{n-1}^{\prime}},
\end{array}
\right.
\end{equation}
%*****************************************************************
where $(\theta_1^{\prime},\dots,\theta_{n-1}^{\prime})$ is the next input vector.
Thus, each vector $(\theta_1,\dots,\theta_{n-1})$ gives rise to a formal 
Jacobi-Perron continued fraction. Whether the fraction is convergent or not, is yet to be 
determined.

Let us introduce the following notation. We let $A^{(0)}=\delta_{ij}$
(the Kronecker delta) and
$
A_i^{(k+n)}=\sum_{j=0}^{n-1} b_i^{(k)}A_i^{(\nu +j)}, \quad b_0^{(k)}=1,
$
where $i=0,\dots,n-1$ and $k=0,1,\dots,\infty$. 
The Jacobi-Perron continued fraction  of the vector 
$(\theta_1,\dots,\theta_{n-1})$  is said to be {\it convergent}, if 
$\theta_i=\lim_{k\to\infty}{A_i^{(k)}\over A_0^{(k)}}$ 
for all $i=1,\dots,n-1$.
Unless $n=2$, convergence of the Jacobi-Perron fractions is a delicate question. 
To the best of our knowledge, there exists no intrinsic necessary 
and sufficient conditions for such a convergence. However,
the Bauer criterion and the Masur-Veech theorem imply
that the Jacobi-Perron fractions  converge for the generic  vectors $(\theta_1,\dots,\theta_{n-1})$.
Namely, let ${\cal F}$ be a measured foliation
on the surface $X$ of genus $g\ge 1$. Recall that the foliation ${\cal F}$ is called uniquely
ergodic if every invariant measure of ${\cal F}$ is a multiple
of the Lebesgue measure. By the Masur-Veech theorem, there exists
a generic subset $V$ in the space of all measured foliations, 
such that each ${\cal F}\in V$ is a uniquely ergodic measured foliation.
We let $\lambda=(\lambda_1,\dots,\lambda_n)$ be the coordinate vector of
the foliation ${\cal F}$. Then the following (Bauer's) criterion 
is true: the Jacobi-Perron continued fraction of $\lambda$  converges 
if and only if $\lambda\in V\subset {\Bbb R}^n$ [Bauer  1996]   \cite{Bau1}.

%**************************************************************************
\subsection{Dimension groups}
%***************************************************************************
 We use ${\Bbb Z},{\Bbb Z}^+,{\Bbb Q}$ and ${\Bbb R}$ for the integers,
positive integers, rationals and reals, respectively and $GL(n,{\Bbb Z})$
for  the group of $n\times n$ matrices with entries in ${\Bbb Z}$
and determinant $\pm 1$.

By an {\it ordered group} we shall mean an abelian group $G$ together
with a subset $P=G^+$ such that $P+P\subseteq P, P\cap (-P)=\{0\}$,
and $P-P=G$. We call $P$ the {\it positive cone} on $G$. We write
$a\le b$ (or $a<b$) if $b-a\in P$ (or $b-a\in P\backslash\{0\}$).

G is said to be a {\it Riesz group} if:

\medskip
(i) $g\in G$ and $ng\ge 0, n\in {\Bbb Z}^+$ implies $g\ge 0$; 

\smallskip
(ii) $u,v\le x,y$ in $G$ implies existence of $w\in G$
such that $u,v\le w\le x,y$.

\smallskip
Given ordered groups $G$ and $H$, we say that a homomorphism
$\varphi:G\to H$ is {\it positive} if $\varphi(G^+)\subseteq H^+$,
and that $\varphi:G\to H$ is an {\it order isomorphism}
if $\varphi(G^+)=H^+$.

A positive homomorphism $f:G\to {\Bbb R}$ is called a {\it state}
if $f(u)=1$, where $u\in G^+$ is an order unit of $G$.
We let $S(G)$ be the {\it state space} of $G$, i.e. the set of states
on $G$ endowed with the natural topology. 

$S(G)$ is a compact convex subset of the vector space $Hom~(G,{\Bbb R})$.
By the Krein-Milman theorem, $S(G)$ is the closed convex hull
of its extreme points, which are called {\it pure states}.

An ordered abelian group is a {\it dimension group}
if it is order isomorphic to $\lim_{m.n\to\infty}
({\Bbb Z}^{r(m)},\varphi_{mn})$, where the ${\Bbb Z}^{r(m)}$
are  simplicially ordered groups (i.e. $({\Bbb Z}^{r(m)})^+
\cong {\Bbb Z}^+\oplus\dots\oplus {\Bbb Z}^+$), and the $\varphi_{mn}$
are positive homomorphisms. The dimension group $G$ is said to be
{\it unimodular} if $r(m)=Const=r$ and $\varphi_{mn}$
are positive isomorphisms of ${\Bbb Z}^r$. In other words, $G$ is
the limit   
%**************************************************************
\begin{equation}\label{eq5}
{\Bbb Z}^r\buildrel\rm
\varphi_0
\over\longrightarrow {\Bbb Z}^r
\buildrel\rm
\varphi_1
\over\longrightarrow
{\Bbb Z}^r\buildrel\rm
\varphi_2
\over\longrightarrow \dots,
\end{equation}
%*********************************************************************
of matrices $\varphi_k\in GL(r,{\Bbb Z}^+)$.

The Riesz groups are dimension groups and vice versa.
The Riesz groups can be viewed as  the abstract dimension groups, while
dimension groups as  a  representation  of the Riesz groups
by the infinite sequences of positive homomorphisms.

\bigskip\noindent
%******************************************************************    
{\sf Acknowledgments.}  I   thank  G.~A.~Elliott  for interesting discussions and   the referee for helpful comments.

%***********************************************************

%**********************************************************

\vskip1cm
\textsc{The Fields Institute for Research in  Mathematical Sciences, Toronto, ON, Canada,  
E-mail:} {\sf igor.v.nikolaev@gmail.com}

\smallskip
{\it Current address: 1505-657 Worcester St.,  Southbridge,  MA 01550,  U.S.A.}

\end{document}